 \documentclass[12pt]{article}
 \usepackage{amsfonts,amssymb,amsmath,mathrsfs}
 \usepackage{color}
 \usepackage[svgnames]{xcolor}

 \newtheorem{thm}{Theorem}[section]
 \newtheorem{lem}[thm]{Lemma}
 \newtheorem{prop}[thm]{Proposition}
 \newtheorem{cor}[thm]{Corollary}

 \def\benumerate{\begin{enumerate}}\def\eenumerate{\end{enumerate}}
 \def\bitemize{\begin{itemize}}\def\eitemize{\end{itemize}}

 \def\beqlb{\begin{eqnarray}}
 \def\eeqlb{\end{eqnarray}}
 \def\beqnn{\begin{eqnarray*}}
 \def\eeqnn{\end{eqnarray*}}

 \def\qqquad{\qquad\qquad}

 \def\proof{\noindent{\it Proof.~~}}
 \def\qed{\hfill$\Box$\medskip}

 \def\<{\langle}\def\>{\rangle}

 \def\mcr{\mathscr}\def\mbb{\mathbb}
 \def\mbf{\mathbf}\def\mrm{\mathrm}\def\mfr{\mathfrak}

 \def\ar{\!\!&}\def\nnm{\nonumber}\def\ccr{\nnm\\}

 \def\d{\mrm{d}}\def\e{\mrm{e}}
 
 \def\supp{\mathrm{supp}}

\begin{document}


\bigskip\bigskip

\centerline{\Large\textbf{Construction of age-structured branching}}

\smallskip

\centerline{\Large\textbf{processes by stochastic equations}}

\bigskip

\centerline{Lina Ji$^{\rm(a)}$ and Zenghu Li$^{\rm(b)}$}

\medskip

\centerline{(a)~Department of Mathematics,}

\centerline{Southern University of Science and Technology,}

\centerline{Shenzhen 518055, People's Republic of China}

\centerline{(b)~Laboratory of Mathematics and Complex Systems,}

\centerline{School of Mathematical Sciences, Beijing Normal University,}

\centerline{Beijing 100875, People's Republic of China}

\centerline{E-mails: \tt jiln@sustech.edu.cn, lizh@bnu.edu.cn}

\bigskip

{\narrower{\narrower

\centerline{\textbf{Abstract}}

\medskip

We give constructions of age-structured branching processes without or with immigration as pathwise unique solutions to stochastic integral equations. A necessary and sufficient condition for the ergodicity of the model with immigration is also given.

\medskip


\par}\par}


\section{Introduction}

Measure-valued branching processes, describing the evolution of population undergoing random reproduction and spatial motion, have been studied by many authors since the pioneer work of Dawson~\cite{D75} and Watanabe~\cite{W68}. We refer to \cite{D93, D94, E00, Li11} for detailed treatments of those processes. A general measure-valued branching process has local and non-local branching mechanisms. Heuristically, the local branching mechanism describes the birth of offspring at the death site of their parent and the non-local one describes that the offspring choose their locations independently according to a distribution as they are born; see, e.g., \cite{DGL02, HLZ15, Li11}, etc. As pointed out by several authors, populations of many species are age-structured. An integer-valued age-dependent branching process was introduced by Bellman and Harris \cite{BH52}; see also \cite{AN72, Har63}. Bose and Kaj \cite{BK00} considered an age-structured branching particle system, where branching events occur during the lifetime and the branching law depends on the age of the parant. A generalized model with spatial motion was introduced in \cite{DGL02}. The reader may also refer to \cite{C94, FHJK20, JK00, JK11} and the references therein for the literature of age-structured models and their properties.

The approach of stochastic equations has played an important role in recent developments of the theory of branching processes. Most of those have been concentrated to the setting of continuous states. In particular, a strong stochastic equation for general continuous-state branching processes with immigration was established by Dawson and Li \cite{DaL06}. A flow of discontinuous continuous-state branching processes was constructed by Bertoin and Le~Gall \cite{BeL06} using weak solutions to a stochastic equation. Their results were extended to general flows in \cite{DaL12, Li14} using strong solutions. For more recent developments in the subject, the reader may refer to \cite{BaM15, Li20, Par16}.

The purpose of this paper is to develop the approach of stochastic equations for the age-structured branching processes. For concreteness, we shall focus to integer-valued processes. To illustrate applications of stochastic equations, we prove a necessary and sufficient condition for the ergodicity of the model with immigration.

Let $\mcr{B}(\mbb{R}_+)$ denote the Borel $\sigma$-algebra on $\mbb{R}_+:= [0,\infty)$. Let $\mfr{M}(\mbb{R}_+)$ denote the set of finite Borel measures on $\mbb{R}_+$ with the weak convergence topology. We define the distance $\rho$ on $\mfr{M}(\mbb{R}_+)$ by
 \beqlb\label{rho}
\rho(\nu_1, \nu_2) = \int_{\mbb{R}_+} e^{-x} |\nu_1(x) - \nu_2(x)| d x,
 \eeqlb
where $\nu_1,\nu_2 \in \mfr{M}(\mbb{R}_+)$. Then $\mfr{M}(\mbb{R}_+)$ is a Polish space whose topology coincides with that given by weak convergence of measures. Let $\mfr{N}(\mbb{R}_+)$ be the subset of $\mfr{M}(\mbb{R}_+)$ consisting of integer-valued measures.

Let $B(\mbb{R}_+)$ be the Banach space of bounded Borel functions on $\mbb{R}_+$ furnished with the supremum norm $\|\cdot\|$. Let $C(\mbb{R}_+)$ be the set of continuous functions in $B(\mbb{R}_+)$ and let $C^1(\mbb{R}_+)$ be the set of functions in $C(\mbb{R}_+)$ with bounded continuous derivatives of the first order. We use the superscript ``+'' to denote the subsets of positive elements of those function spaces, e.g., $B(\mbb{R}_+)^+$, $C(\mbb{R}_+)^+$, etc. For any $f \in B(\mbb{R}_+)$ and $\nu \in \mfr{M}(\mbb{R}_+)$ write $\<\nu,f\> = \int_{\mbb{R}_+} fd\nu$.

Throughout the paper, we let the letter $K$ with or without subscripts to denote constants whose exact value is unimportant and may change from line to line. In the integrals, we make the convention that, for $a \le b \in \mbb{R},$
 \beqnn
\int_a^b = \int_{(a, b]} \quad \text{and}\quad \int_a^\infty = \int_{(a, \infty)}.
 \eeqnn

The rest of this paper is organized as follows. In Section 2, the model of age-structured branching process is constructed as the unique strong solution of a stochastic integral equation. The characterizations of some interesting probabilities related to the population evolution are given. Similar results for the age-structured system with immigration are presented in Section 3, where the ergodicity of the model is also studied.

\section{An age-structured branching process}

 \setcounter{equation}{0}

In this section, we give formulations of the age-structured branching process in terms of some stochastic integral equations. We shall prove the existence and uniqueness of the solutions to the equations. Let $\alpha\in B(\mbb{R}_+)^{+}$ be bounded away from zero and let $g=g(x,z)$ be a positive Borel function on $\mbb{R}_+\times [0,1]$ of the form:
 \beqnn
g(x,z) = \sum_{i=0}^\infty p(x, i)z^i,
 \eeqnn
where $p(x, i)\ge 0$ and $\sum_{i=0}^\infty p(x, i)=1$. Then $g(x,\cdot)$ is a probability generating function for any fixed $x\in \mbb{R}_+$. We assume that
 \beqlb\label{2.1zz}
\|g'(\cdot,1-)\| = \sup_{x\ge 0}\sum_{i=1}^\infty p(x,i)i< \infty.
 \eeqlb
A branching particle system is characterized by the following properties:
 \benumerate

\item[(2.A)] The ages of the particles increase at the unit speed, i.e., they move as the Markov process $\xi= (\xi_t)_{t\ge 0}$ in $\mbb{R}_+$ defined by $\xi_t= \xi_0+t$.

\item[(2.B)] For a particle which is alive at time $r\ge 0$ with age $x\ge 0$, the conditional probability of survival in the time interval $[r,t)$ is $\exp\{-\int_0^{t-r} \alpha(x+s)d s\}$.

\item[(2.C)] When a particle dies at age $x\ge 0$, it gives birth to a random number of offspring with age zero according to the probability law given by the generating function $g(x,\cdot)$.

 \eenumerate
We assume that the lifetimes and the offspring productions of different particles are independent. Let $X_t(B)$ denote the number of particles with ages belonging the Borel set $B\subset \mbb{R}_+$ that are alive at time $t\ge 0$. If we assume $X_0(\mbb{R}_+)< \infty$, then $\{X_t: t\ge 0\}$ is a Markov process with state space $\mfr{N}(\mbb{R}_+)$. We refer to Li~\cite[Section~4.3]{Li11} for the formulation of general branching particle systems.

Let $\sigma\in \mfr{N}(\mbb{R}_+)$ and let $\{X_t^\sigma: t\ge 0\}$ be the above system with initial value $X^\sigma_0 = \sigma$. Suppose that the process is defined on a probability space $(\Omega, \mcr{G}, \mbf{P})$. The above properties imply that
 \beqlb\label{2.1}
\mbf{E}\exp\{-\<X_t^\sigma, f\>\} = \exp\{-\<\sigma, u_tf\>\}, \qquad f\in B(\mbb{R}_+)^+,
 \eeqlb
where
 \beqnn
u_tf(x) = -\log \mbf{E}\exp\{-\<X_t^{\delta_x}, f\>\}.
 \eeqnn
From properties (2.A), (2.B) and (2.C) we derive as in Li~\cite[p.89]{Li11} the following renewal equation:
{\small \beqlb\label{2.2}
\e^{-u_tf(x)} = \e^{-f(x+t)-\int_0^t \alpha(x+s)d s} + \int_0^t \e^{-\int_0^s\alpha(x+r)d r}\alpha(x+s) g(x+s,\e^{-u_{t-s}f(0)})d s.
 \eeqlb}
By Li~\cite[Proposition~2.9]{Li11}, the above equation implies
 \beqlb\label{2.3}
\e^{-u_tf(x)} = \e^{-f(x+t)} + \int_0^t \alpha(x+s)\big[g(x+s,\e^{-u_{t-s}f(0)}) - \e^{-u_{t-s}f(x+s)}\big]d s.
 \eeqlb
The uniqueness of the solution to \eqref{2.2} and \eqref{2.3} follows by Gronwall's inequality.

We call any Markov process $(X_t: t\ge 0)$ with state space $\mfr{N}(\mbb{R}_+)$ an \textit{$(\alpha,g)$-age-structured branching process} if it has transition semigroup $(Q_t)_{t\ge 0}$ defined by
 \beqlb\label{2.4}
\int_{\mfr{N}(\mbb{R}_+)}\e^{-\<\nu, f\>}Q_t(\sigma,d\nu) = \exp\{-\<\sigma, u_tf\>\}, \qquad f\in B(\mbb{R}_+)^+,
 \eeqlb
where $u_tf(x)$ is the unique solution to \eqref{2.3}.

\begin{prop}\label{t2.1a}
For any $t\ge 0$ and $\sigma\in \mfr{N}(\mbb{R}_+)$ we have
 \beqlb\label{2.5}
\int_{\mfr{N}(\mbb{R}_+)}\<\nu,f\> Q_t(\sigma,d\nu)
 =
\<\sigma, \pi_tf\>, \qquad f\in B(\mbb{R}_+),
 \eeqlb
where $(\pi_t)_{t\ge 0}$ is the semigroup of bounded kernels on $\mbb{R}_+$ defined by
 \beqlb\label{2.6}
\pi_tf(x) = f(x+t) + \int_0^t \alpha(x+s)[g'(x+s,1-)\pi_{t-s}f(x+s) - \pi_{t-s}f(0)] d s.
 \eeqlb
\end{prop}

\proof The existence and uniqueness of the locally bounded solution to \eqref{2.6} follows by a general result; see, e.g., Lemma~2.17 of Li (2011). For $f\in B(\mbb{R}_+)^+$ we can use \eqref{2.3} to see the unique solution of \eqref{2.6} is given by $\pi_tf(x)= (d/d\theta)u_t(\theta f)(x)|_{\theta=0}$. By differentiating both sides of \eqref{2.4} we get \eqref{2.5}. The extension to $f\in B(\mbb{R}_+)$ is immediate by linearity. \qed

\begin{prop}\label{t2.1b}
For any $t\ge 0$ and $x\ge 0$ we have
 \beqlb\label{2.5b}
 u_tf(x) \ge (1-\e^{-f(x+t)})e^{-c_1t}, \qquad f\in B(\mbb{R}_+)^+,
 \eeqlb
where $c_1= \sup_{y\ge 0}\alpha(y)$. \end{prop}

\proof Let $U_tf(x)= 1 - e^{-u_tf(x)}$. From \eqref{2.3} we obtain
 \beqnn
U_tf(x) = 1-\e^{-f(x+t)} + \int_0^t \alpha(x+s)\big[h(s,t,x) - U_{t-s}f(x+s)\big]d s,
 \eeqnn
where
 \beqnn
h(s,t,x) = 1-g(x+s,\e^{-u_{t-s}f(0)})\ge 0.
 \eeqnn
By comparison theorem we have $U_tf(x)\ge v_tf(x)$, where $(t,x)\mapsto v_tf(x)$ solves
 \beqnn
v_tf(x) = 1-\e^{-f(x+t)} - \int_0^t \alpha(x+s)v_{t-s}f(x+s) d s.
 \eeqnn
The unique locally bounded solution to the above equation is given by
 \beqnn
v_tf(x) = (1-\e^{-f(x+t)})\exp\bigg\{-\int_0^t \alpha(x+s) d s\bigg\}.
 \eeqnn
Then we have the estimate \eqref{2.5b}. \qed

Let $D(\mbb{R}_+)$ be the set of bounded positive right-continuous increasing functions $f$ on $\mbb{R}$ satisfying $f(x)= 0$ for $x< 0$. We identify $\mu\in \mfr{M}(\mbb{R}_+)$ with its distribution function $\mu\in D(\mbb{R}_+)$ defined by $\mu(x)= \mu[0,x]$ for $x\ge 0$. For $\mu\in D(\mbb{R}_+)$ and $x\in \mbb{R}$ let $\mu^{-1}(x)= \inf\{y\ge 0: \mu(y)> x\}$ and $\tilde{\mu}^{-1}(x)= \mu^{-1}(\mu(\infty)x)$, where $\mu(\infty)= \lim_{x\to \infty} \mu(x)= \mu(\mbb{R}_+)$. For $\mu\in D(\mbb{R}_+)$ and $\alpha\in B(\mbb{R}_+)^+$ we define ${_\alpha\mu}\in D(\mbb{R}_+)$ by ${_\alpha\mu}(x)= \int_{[0,x]}\alpha(y)\mu(dy)$ for $x\ge 0$.

Suppose that $(\Omega, \mcr{F}, \mcr{F}_t, \mbf{P})$ is a filtered probability space satisfying the usual hypotheses. Let $M(d t, d u, d y, d z, d v)$ be an $(\mcr{F}_t)$-Poisson random measure on $(0,\infty)^2\times (0,1]\times \mbb{N}\times (0,1]$ with intensity $dtdudy\pi(dz)dv$, where $\pi(dz)$ denotes the counting measure on $\mbb{N}$. Given an $\mcr{F}_0$-measuable random function $X_0\in D(\mbb{R}_+)$, we consider the following stochastic integral equation:
{\small \beqlb\label{2.4added01}
X_t(x) \ar=\ar X_0(x-t) + \int_0^t\int_0^{\<X_{s-},\alpha\>} \int_0^1\int_{\mbb{N}}\int_0^{p({_\alpha \tilde{X}}_{s-}^{-1}(y), z)} z1_{\{t-s\le x\}} M(ds,du,dy,dz,d v) \cr
 \ar\ar
- \int_0^t\int_0^{\<X_{s-},\alpha\>}\int_0^1\int_{\mbb{N}} \int_0^{p({_\alpha\tilde{X}}_{s-}^{-1}(y), z)} 1_{\{{_\alpha \tilde{X}}_{s-}^{-1}(y) + t-s\le x\}} M(ds,du,dy,dz,dv). \qquad
 \eeqlb}
Heuristically, the left-hand side $X_t(x)$ is the number of individuals at time $t$ with ages less than $x$. On the right-hand side, the first term $X_0(x-t)$ counts the number of individuals having ages less than $x-t$ at time $0$ and having ages less than $x$ at time $t$. A death of the population occurs at time $s\in [0,t]$ at rate $\<X_{s-}, \alpha\>ds$. In that case, the age of the dying individual is distributed according to the probability measure $ \<X_{s-}, \alpha\>^{-1}{_\alpha X}_{s-}$ and is realized as
 \beqnn
{_\alpha \tilde{X}}_{s-}^{-1}(y)= {_\alpha X}_{s-}^{-1}(\<X_{s-}, \alpha\>y)
 =
\inf\{z\ge 0: \<X_{s-}, \alpha\>^{-1}{_\alpha X}_{s-}(z)>y\},
 \eeqnn
where $y$ is selected from $[0,1]$ according to the uniform distribution by the Poisson random measure. The number of offspring of the individual takes the value $z\in \mbb{N}$ at ratio $p({_\alpha \tilde{X}}_{s-}^{-1}(y), z)$ and contributes to the number $X_t(x)$ if and only if $t-s\le x$, which is recorded by the second term. The death of the individual makes effect on $X_t(x)$ if and only if ${_\alpha \tilde{X}}_{s-}^{-1}(y) + t-s\le x$, which is recorded by the third term.

Let $\zeta_a(x) = 1_{\{a\le x\}}$ for $a, x\in \mbb{R}$. Given a function $f$ on $\mbb{R}$ define $f\circ\theta_t(x)= f(x + t)$ for $x, t\in \mbb{R}$. Then we may rewrite \eqref{2.4added01} equivalently into
{\small \beqlb\label{2.4added02}
X_t(x) \ar=\ar X_0\circ\theta_{-t}(x) + \int_0^t\int_0^{\<X_{s-},\alpha\>}\int_0^1\int_{\mbb{N}}\int_0^{p({_\alpha \tilde{X}}_{s-}^{-1}(y), z)} \Big[z \zeta_0\circ\theta_{s-t}(x) \cr
 \ar\ar\qqquad\qqquad\qquad
-\, \zeta_{{_\alpha \tilde{X}}_{s-}^{-1}(y)}\circ\theta_{s-t}(x)\Big] M(ds,du,dy,dz,dv).
 \eeqlb}
A pathwise unique solution to \eqref{2.4added02} is constructed as follows. Let $\tau_0 = 0$. Given $\tau_{k-1}\ge 0$ and $X_{\tau_{k-1}}\in D(\mbb{R}_+)$, we first define
 \beqnn
\tau_k = \tau_{k-1} + \inf\{t>0: M((\tau_{k-1},\tau_{k-1}+t]\times (0,\<X_{\tau_{k-1}}, \alpha\>]\times A_k)> 0\},
 \eeqnn
where $A_k= \{(y,z,v): y\in (0,1], z\in\mbb{N},0< v\le p({_\alpha \tilde{X}}_{\tau_{k-1}}^{-1}(y),z)\}$, and
 \beqlb\label{2.4added02a}
X_t(x)= X_{\tau_{k-1}}\circ\theta_{\tau_{k-1}-t}(x), \qquad \tau_{k-1}< t< \tau_k.
 \eeqlb
Then we define
 \beqlb\label{2.4added02b}
X_{\tau_k}(x)= X_{\tau_k-}(x) + z_k\zeta_0(x) - \zeta_{{_\alpha\tilde{X}}_{\tau_k-}^{-1}(y_k)}(x),
 \eeqlb
where $X_{\tau_k-}(x)= X_{\tau_{k-1}}\circ\theta_{\tau_{k-1}-\tau_k}(x)$ and $(u_k,y_k,z_k,v_k)\in (0,\infty)\times (0,1]\times \mbb{N}\times (0,1]$ is the point so that $(\tau_k, u_k,y_k,z_k,v_k)\in \supp(M)$. Since ${_\alpha\tilde{X}}_{\tau_k-}^{-1}(y_k)\in \supp(X_{\tau_k-})$, we have $X_{\tau_k}\in D(\mbb{R}_+)$. The expression \eqref{2.4added02b} means that at time $\tau_k$ an individual at ${_\alpha\tilde{X}}_{\tau_k-}^{-1}(y_k)$ dies and gives birth to $z_k$ offspring with age $0\in \mbb{R}_+$. It is easy to see that $\{X_t: t \ge 0\}$ is the pathwise unique solution to \eqref{2.4added01} or \eqref{2.4added02} with lifetime $\tau:= \lim_{k\to \infty}\tau_k$. More precisely, the equations hold a.s.\ with $t$ replaced by $t\land \tau_k$ for every $k\ge 1$. Let $n(t)= \sup\{k\ge 0: \tau_k\le t\}$ for $t\ge 0$. Recall that $\beta:= \|\alpha g'(\cdot,1-)\|< \infty$.

\begin{lem}\label{t2.2a2s1}
Suppose that $\mbf{E}[X_0(\infty)]< \infty$. Then for any $k\ge 1$ we have
 \beqlb\label{2.15a1aa}
\mbf{E}\Big[\sup_{0\le s\le t\land\tau_k}X_s(\infty)\Big]
 \le
\mbf{E}[X_0(\infty)]\e^{\beta t}, \qquad t\ge 0.
 \eeqlb
\end{lem}

\proof Recall that $X_t(\infty)= \lim_{x\to \infty} X_t(x)= X_t(\mbb{R}_+)$. In view of \eqref{2.4added02}, we have
{\small \beqlb\label{2.4added02axx}
X_t(\infty)= X_0(\infty) + \int_0^t\int_0^{\<X_{s-},\alpha\>} \int_0^1\int_{\mbb{N}}\int_0^{p({_\alpha \tilde{X}}_{s-}^{-1}(y),z)} (z-1) M(ds,du,dy,dz,dv).
 \eeqlb}
It follows that
{\small \beqnn
\mbf{E}\Big[\sup_{0\le s\le t\land\tau_k}X_s(\infty)\Big]
 \ar\le\ar
\mbf{E}[X_0(\infty)] + \sum_{z\in \mbb{N}} \mbf{E}\bigg[\int_0^{t\land\tau_k} \<X_{s-},\alpha\>d s \int_0^1 p({_\alpha \tilde{X}}_{s-}^{-1}(y),z) z d y\bigg] \cr
 \ar=\ar
\mbf{E}[X_0(\infty)] + \sum_{z\in \mbb{N}} \mbf{E}\bigg[\int_0^{t\land\tau_k} d s \int_{\mbb{R}_+} \alpha(y)p(y,z)z dX_{s-}(y)\bigg] \cr
 \ar=\ar
\mbf{E}[X_0(\infty)] + \mbf{E}\bigg[\int_0^{t\land\tau_k} d s \int_{\mbb{R}_+} \alpha(y)g'(y,1-) dX_{s-}(y)\bigg] \cr
 \ar\le\ar
\mbf{E}[X_0(\infty)] + \beta\mbf{E}\bigg[\int_0^{t\land\tau_k} X_{s-}(\infty)d s\bigg].
 \eeqnn}
Then {\small $\mbf{E}[\sup_{0\le s\le t\land\tau_k}X_s(\infty)]$} is locally bounded in $t\ge 0$ and
{\small \beqnn
\mbf{E}\Big[\sup_{0\le s\le t\land\tau_k}X_s(\infty)\Big]
 \ar\le\ar
\mbf{E}[X_0(\infty)] + \beta\mbf{E}\bigg[\int_0^{t\land\tau_k} X_s(\infty)d s\bigg] \cr
 \ar\le\ar
\mbf{E}[X_0(\infty)] + \beta\int_0^t \mbf{E}[X_{s\land\tau_k}(\infty)]d s \cr
 \ar\le\ar
\mbf{E}[X_0(\infty)] + \beta\int_0^t \Big[\sup_{0\le r\le s\land\tau_k}X_r(\infty)\Big]d s.
 \eeqnn}
By Gronwall's inequality we obtain \eqref{2.15a1aa}. \qed

\begin{lem}\label{t2.2a2s2}
Suppose that $\mbf{E}[X_0(\infty)]< \infty$. Then we have $\mbf{P}\{\tau= \infty\}= 1$ and
 \beqlb\label{2.15a1bb}
\mbf{E}[n(t)]\le \|\alpha\|\mbf{E}[X_0(\infty)] \int_0^t \e^{\beta s} d s, \qquad t\ge 0.
 \eeqlb
\end{lem}

\proof By \eqref{2.4added02} and monotone convergence we have
 \beqnn
\mbf{E}[n(t)]\ar=\ar \lim_{k\to \infty}\mbf{E}\bigg[\int_0^{t\land\tau_k}\int_0^{\<X_{s-},\alpha\>} \int_0^1\int_{\mbb{N}}\int_0^{p({_\alpha \tilde{X}}_{s-}^{-1}(y),z)} M(ds,du,dy,dz,dv)\bigg] \cr
 \ar=\ar
\lim_{k\to \infty}\mbf{E}\bigg[\sum_{z\in \mbb{N}} \int_0^{t\land\tau_k} \<X_{s-},\alpha\> d s \int_0^1 p({_\alpha \tilde{X}}_{s-}^{-1}(y),z) dy\bigg] \cr
 \ar=\ar
\lim_{k\to \infty}\mbf{E}\bigg[\sum_{z\in \mbb{N}} \int_0^{t\land\tau_k} \alpha(y) d s \int_{\mbb{R}_+} p(y,z) d\tilde{X}_{s-}(y)\bigg] \cr
 \ar\le\ar
\lim_{k\to \infty}\|\alpha\|\mbf{E}\bigg[\int_0^{t\land\tau_k} X_s(\infty) d s\bigg]
 \le
\lim_{k\to \infty}\|\alpha\|\int_0^t \mbf{E}[X_{s\land\tau_k}(\infty)] d s.
 \eeqnn
Then \eqref{2.15a1bb} follows by \eqref{2.15a1aa}. In particular, we have $\mbf{P}\{\tau> t\}= \mbf{P}\{n(t)< \infty\}= 1$ for every $t\ge 0$. That implies $\mbf{P}\{\tau= \infty\}= 1$. \qed

\begin{prop}\label{t2.2a3}
Suppose that $\mbf{E}[X_0(\infty)]< \infty$. Then we have
 \beqnn
\mbf{E}\Big[\sup_{0\le s\le t}X_s(\infty)\Big]
 \le
\mbf{E}[X_0(\infty)]\e^{\beta t}, \qquad t\ge 0.
 \eeqnn
\end{prop}

\proof Since $\mbf{P}\{\tau= \infty\}= 1$ by Lemma~\ref{t2.2a2s2}, we obtain the result from \eqref{2.15a1aa} by using monotone convergence. \qed

By Lemma~\ref{t2.2a2s2} the solution of \eqref{2.4added01} or \eqref{2.4added02} has infinite lifetime and determines a measure-valued strong Markov process $\{X_t: t\ge 0\}$. A characterization of this process is given by the following:

\begin{prop}\label{t2.2a1}
For any $t\ge 0$ and $f \in B(\mbb{R}_+)$ we have
 \beqlb\label{2.15}
\<X_t,f\> \ar=\ar \<X_0,f\circ \theta_t\> + \int_0^t\int_0^{\<X_{s-},\alpha\>}\int_0^1 \int_{\mbb{N}} \int_0^{p({_\alpha \tilde{X}}_{s-}^{-1}(y),z)} \Big[z f\circ\theta_{t-s}(0) \cr
 \ar\ar\qqquad\qqquad
-\, f\circ\theta_{t-s}({_\alpha \tilde{X}}_{s-}^{-1}(y))\Big] M(ds,du,dy,dz,dv). \quad
 \eeqlb
\end{prop}

\proof Let $C^1_0(\mbb{R}_+)$ denote the subspace of $C^1(\mbb{R}_+)$ consisting of functions vanishing at infinity. For any fixed integer $n\ge 1,$ let $x_i = in/2^n$ with $i = 0, 1, \cdots, 2^n.$ By \eqref{2.4added02} it holds almost surely for any $f \in C^1(\mbb{R}_+)$ that
 \beqnn
\ar\ar\sum_{i = 1}^{2^n}f'(x_i)X_t(x_i) - \sum_{i = 1}^{2^n}f'(x_i)X_0\circ\theta_{- t}(x_i) \cr
 \ar\ar\qquad
= \int_0^t\int_0^{\<X_{s-}, \alpha\>}\int_0^1\int_{\mbb{N}}\int_0^{p({_\alpha \tilde{X}}_{s-}^{-1}(y), z)}\sum_{i = 1}^{2^n}f'(x_i)\Big[z \zeta_0\circ\theta_{s-t}(x_i) \cr
 \ar\ar\qqquad\qqquad\qqquad\qquad
-\, \zeta_{Z_{s-}(y)}\circ\theta_{s- t}(x_i)\Big] M(ds,du,dy,dz,dv).
 \eeqnn
Then we multiply the above equation by $2^{-n}$ and let $n\to \infty$ to see, almost surely,
 \beqlb\label{2.14}
\ar\ar\int_0^\infty f'(x)X_t(x)d x - \int_0^\infty f'(x)X_0\circ\theta_{-t}(x)d x \cr
 \ar\ar\qquad
= \int_0^t\int_0^{\<X_{s-}, \alpha\>}\int_0^1\int_{\mbb{N}}\int_0^{p({_\alpha \tilde{X}}_{s-}^{-1}(y), z)} \bigg\{\int_0^\infty f'(x) \Big[z\zeta_0\circ\theta_{s-t}(x) \cr
 \ar\ar\qqquad\qqquad\qqquad
-\, \zeta_{{_\alpha \tilde{X}}_{s-}^{-1}(y)}\circ\theta_{s-t}(x)\Big]\d x\bigg\} M(ds,du,dy,dz,dv) \cr
 \ar\ar\qquad
= \int_0^t\int_0^{\<X_{s-}, \alpha\>}\int_0^1\int_{\mbb{N}}\int_0^{p({_\alpha \tilde{X}}_{s-}^{-1}(y), z)} \Big[z f\circ\theta_{t-s}(0) \cr
 \ar\ar\qqquad\qqquad\qqquad
-\, f\circ\theta_{t-s}({_\alpha \tilde{X}}_{s-}^{-1}(y))\Big] M(ds,du,dy,dz,dv). \quad
 \eeqlb
By integration by parts we have
 \beqlb\label{2.13}
\<X_t, f\> = - \int_0^\infty f'(x)X_t(x)d x.
 \eeqlb
From \eqref{2.14} and \eqref{2.13} we see that \eqref{2.15} holds for any $f \in C^1(\mbb{R}_+)$. Then the relation also holds for any $f\in B(\mbb{R}_+)$ by a monotone class argument. \qed

\begin{prop}\label{t2.2a5}
For any $t\ge 0$ and $f\in C^1(\mbb{R}_+)$ we have
{\small \beqlb\label{2.15a1}
\<X_t, f\> \ar=\ar \<X_0, f\> + \int_0^t \<X_{s-}, f'\> d s + \int_0^t \int_0^{\<X_{s-},\alpha\>}\int_0^1\int_{\mbb{N}}\int_0^{p({_\alpha \tilde{X}}_{s-}^{-1}(y), z)}\Big[z f(0) \cr
 \ar\ar\qqquad\qqquad\qqquad
-\, f({_\alpha \tilde{X}}_{s-}^{-1}(y))\Big] M(ds,du,dy,dz,dv). \qquad
 \eeqlb}
\end{prop}

\proof For $n\ge 1$ we consider a partition $\Delta_n = \{0 = t_0 < t_1 < \cdots < t_n = t\}$ of $[0, t].$ Notice that $df\circ \theta_t(x)/dt = f'(x+t).$ By \eqref{2.15} we have
 \beqnn
\<X_t, f\> \ar=\ar \<X_0, f\> + \sum_{i = 1}^n \left[\<X_{t_i}, f\> - \<X_{t_{i - 1}}, f\circ\theta_{t_{i} - t_{i - 1}}\>\right]\cr
\ar\ar + \sum_{i = 1}^n \left[\<X_{t_{i - 1}}, f\circ\theta_{t_{i} - t_{i - 1}}\> - \<X_{t_{i - 1}},f\>\right]\cr
\ar=\ar \<X_0, f\> + \sum_{i = 1}^n \int_{t_{i - 1}}^{t_i} \int_0^{\<X_{s-},\alpha\>}\int_0^1\int_{\mbb{N}}\int_0^{p({_\alpha \tilde{X}}_{s-}^{-1}(y), z)} \Big[z f\circ\theta_{t_i-s}(0) \cr
 \ar\ar
-\, f\circ\theta_{t_i - s}({_\alpha \tilde{X}}_{s-}^{-1}(y))\Big] M(ds,du,dy,dz,dv)\cr
\ar\ar + \sum_{i = 1}^n \int_{t_{i - 1}}^{t_i}\<X_{t_{i - 1}}, f'\circ\theta_{s - t_{i - 1}}\>d s.
 \eeqnn
By letting $|\Delta_n|:= \max_{1\le i\le n} (t_n-t_{n-1})\to 0$ and using the right continuity of $s\to X_s$ and the continuity of $s\to f\circ\theta_s$ we obtain \eqref{2.15a1}. \qed

 \begin{prop}\label{t2.2a4}
For $f, G \in C^1(\mbb{R}_+)$ let $G_f(\mu)= G(\<\mu,f\>)$ and
 \beqnn
LG_f(\mu)\ar=\ar \<\mu,f'\>G'(\<\mu,f\>) - \int_{\mbb{R}_+} \alpha(y) \sum_{z \in \mbb{N}}p(y,z)\big[G(\<\mu, f\>) \ccr
 \ar\ar\qqquad
-\, G\big(\<\mu,f\>+zf(0)-f(y)\big)\big] \mu(d y).
 \eeqnn
Then we have
 \beqlb\label{2.15a4}
G_f(X_t) = G_f(X_0) + \int_0^t LG_f(X_s)d s + \mbox{mart.}
 \eeqlb
\end{prop}

\proof Let $\tilde{M}$ denote the compensated measure of $M$. By Proposition~\ref{t2.2a5} and It\^{o}'s formula, we have
 \beqnn
G(\<X_t,f\>)\ar=\ar G(\<X_0,f\>) + \int_0^t G'(\<X_s,f\>) \<X_s, f'\>d s \cr
 \ar\ar
- \int_0^t \int_0^{\<X_{s-},\alpha\>}\int_0^1\int_{\mbb{N}}\int_0^{p({_\alpha \tilde{X}}_{s-}^{-1}(y),z)} \Big[G(\<X_{s-},f\>) \ccr
 \ar\ar\qquad
-\,G(\<X_{s-},f\> + zf(0) - f({_\alpha \tilde{X}}_{s-}^{-1}(y)))\Big] M(ds,du,dy,dz,dv) \cr
\ar=\ar
G(\<X_0,f\>) + \int_0^t G'( \<X_s, f\>) \<X_s, f'\>d s - M_t^G(f) \cr
 \ar\ar
- \int_0^tds \int_{\mbb{R}_+} \alpha(y)\sum_{z\in \mbb{N}}p(y,z) \Big[G(\<X_{s-},f\>) \ccr
 \ar\ar\qqquad\qqquad
-\,G\big(\<X_{s-},f\> + zf(0) - f(y)\big)\Big] X_{s-}(d y),
 \eeqnn
where
 \beqnn
M_t^G(f)\ar=\ar \int_0^t \int_0^{\<X_{s-},\alpha\>}\int_0^1\int_{\mbb{N}}\int_0^{p({_\alpha \tilde{X}}_{s-}^{-1}(y), z)}\Big[G\big(\<X_{s-}, f\>\big) \ccr
 \ar\ar\qquad
-\,G\big(\<X_{s-}, f\> + z f(0) - f({_\alpha \tilde{X}}_{s-}^{-1}(y))\big)\Big] \tilde{M}(ds,du,dy,dz,dv).
 \eeqnn
By Proposition~\ref{t2.2a3} one can check that $\{M_t^G(f): t\ge 0\}$ is a martingale.
 \qed

\begin{thm}\label{t2.2a}
The measure-valued process $\{X_t: t \ge 0\}$ defined by \eqref{2.4added01} is an $(\alpha, g)$-age-structured branching process.
\end{thm}

\proof
Let $G \in C^1(\mbb{R}_+)$ and let $t\to f_t$ be a mapping from $[0, T]$ to $C^1(\mbb{R}_+)^+$ such that $t\to f_t$ is continuously differentiable and $t\to f_t'$ is continuous by the supremum norm. For $0 \le t \le T$ and $k \ge 1$ we have
 \beqnn
G(\<X_t,f_t\>) \ar=\ar G(\<X_0,f_0\>) + \sum_{j = 0}^\infty \Big[G(\<X_{t\wedge (j + 1)/k},f_{t\wedge j/k}\>) - G(\<X_{t\wedge j/k},f_{t\wedge j/k}\>)\Big] \cr
\ar\ar + \sum_{j = 0}^\infty \Big[G(\<X_{t\wedge (j + 1)/k},f_{t\wedge (j + 1)/k}\>) - G(\<X_{t\wedge (j + 1)/k},f_{t\wedge j/k}\>)\Big],
 \eeqnn
where the summations only consist of finitely many non-trivial terms. Based on Proposition~\ref{t2.2a4}, we obtain
 {\small \beqnn
G(\<X_t,f_t\>) \ar=\ar G(\<X_0,f_0\>) + \sum_{j = 0}^\infty \int_{t\wedge j/k}^{t\wedge (j + 1)/k} \Big\{ G'(\<X_s,f_{t\wedge j/k}\>)\<X_s,f_{t\wedge j/k}'\> + \int_{\mbb{R}_+}\alpha(y)X_s(d y) \cr
\ar\ar \sum_{z \in \mbb{N}}p(y,z)\Big[G(\<X_s,f_{t\wedge j/k}\> + zf_{t\wedge j/k}(0) - f_{t\wedge j/k}(y)) - G(\<X_s,f_{t\wedge j/k}\>)\Big]\Big\}d s \cr
\ar\ar + \sum_{j = 0}^\infty \int_{t\wedge j/k}^{t\wedge (j + 1)/k}G' (\<X_{t\wedge (j + 1)/k},f_s\>)\<X_{t\wedge (j + 1)/k}, d f_s/d s\>ds + M_k(t),
 \eeqnn\!\!}
where $\{M_k(t)\}$ is a martingale. Recall that $\{X_t\}$ is a c\'{a}dl\'{a}g process, letting $k\to \infty$ in the above equation we have
 \beqnn
G(\<X_t,f_t\>) \ar=\ar G(\<X_0,f_0\>) + \int_0^t \Big\{ G'(\<X_s,f_s\>)\<X_s,f_s'\> + \int_{\mbb{R}_+}\alpha(y)X_s(d y) \cr
\ar\ar \sum_{z \in \mbb{N}}p(y,z)\Big[G(\<X_s,f_s\> + zf_s(0) - f_s(y)) - G(\<X_s,f_s\>)\Big]\Big\}d s \cr
\ar\ar + \int_0^t G'(\<X_s,f_s\>)\<X_s, d f_s/d s\> d s + M(t),
 \eeqnn
where $\{M(t)\}$ is a martingale. For any $f \in C^1(\mbb{R}_+)^+$ we may apply the above to $G(z) = e^{-z}$ and $f_t = u_{T -t}f$, where $u_tf$ is the unique solution to \eqref{2.3}. Then one may see that $t\to \exp\{- \<X_t, u_{T-t}f\>\}$ is a martingale, so $\{X_t: t \ge 0\}$ has transition semigroup $(Q_t)_{t\ge 0}$ defined by \eqref{2.4}. It follows that $\{X_t: t \ge 0\}$ is an $(\alpha, g)$-age-structured branching process. \qed

\begin{prop}\label{t2.2a4x}
Suppose that $\mbf{E}[\<X_0,1\>]< \infty$. Then for any $t\ge 0$ and $f\in B(\mbb{R}_+)^+$ we have
 \beqlb\label{2.15a2x}
\mbf{E}[\<X_t,1\>]
 \le
\mbf{E}[\<X_0,1\>]\e^{c_0t},
 \eeqlb
where
 \beqlb\label{2.15a2zz}
c_0= \sup_{y\ge 0}\alpha(y)\big[g'(y,1-)-1\big].
 \eeqlb
\end{prop}

\proof The equality in \eqref{2.15a2x} follows from Proposition~\ref{t2.1a} and the Markov property. Since $X_t(\infty)= \<X_t,1\>$, by taking expectations in \eqref{2.4added02axx} we obtain
 \beqnn
\mbf{E}[\<X_t,1\>]
 \ar=\ar
\mbf{E}[\<X_0,1\>] + \sum_{z\in \mbb{N}} \mbf{E}\bigg[\int_0^t \<X_{s-},\alpha\>d s \int_0^1 p({_\alpha \tilde{X}}_{s-}^{-1}(y),z) (z-1) d y\bigg] \cr
 \ar=\ar
\mbf{E}[\<X_0,1\>] + \sum_{z\in \mbb{N}} \mbf{E}\bigg[\int_0^t d s \int_{\mbb{R}_+} \alpha(y)p(y,z) (z-1) dX_s(y)\bigg] \cr
 \ar=\ar
\mbf{E}[\<X_0,1\>] + \mbf{E}\bigg[\int_0^t d s \int_{\mbb{R}_+} \alpha(y)[g'(y,1-)-1] dX_s(y)\bigg].
 \eeqnn
Then we get \eqref{2.15a2x} by a comparison result. \qed

\begin{cor}\label{t2.1c}
We have
 \beqlb\label{2.5c}
\|\pi_tf\|\le e^{c_0t}\|f\|, \qquad t\ge 0, f\in B(\mbb{R}_+).
 \eeqlb
\end{cor}

\proof By \eqref{2.5} and the Markov property we have $\mbf{E}[\<X_t,f\>]= \mbf{E}[\<X_0,\pi_tf\>]$. But from the estimate given in Proposition~\ref{t2.2a4x} it follows that
 \beqnn
|\mbf{E}[\<X_t,f\>]|
 \le
\|f\|\mbf{E}[\<X_t,1\>]
 \le
\|f\|\mbf{E}[\<X_0,1\>]\e^{c_0t}.
 \eeqnn
By letting $X_0=\delta_x$ for $x\ge 0$ we see $|\pi_tf(x)|\le e^{c_0t}\|f\|$, which implies \eqref{2.5c}. \qed

\section{The branching process with immigration}

 \setcounter{equation}{0}

In this section we introduce a generalization of the model discussed in Section~2, which is the age-structured branching process with immigration. Let $\psi$ be a functional on $B(\mbb{R}_+)^+$ given by
 \beqnn
\psi(f) = \int_{\mfr{N}(\mbb{R}_+)^\circ} (1 - e^{-\<\nu, f\>})L (d\nu), \qquad f\in B(\mbb{R}_+)^+,
 \eeqnn
where $L(d\nu)$ is a finite positive measure on $\mfr{N}(\mbb{R}_+)^\circ := \mfr{N}(\mbb{R}_+)\setminus 0,$ and $0$ denotes the null measure.

A Markov process $Y= (Y_t: t\ge 0)$ with state space $\mfr{N}(\mbb{R}_+)$ is called an \textit{$(\alpha,g, \psi)$-age-structured branching process with immigration} if it has transition semigroup $(P_t)_{t\ge 0}$ defined by:
 \beqlb\label{1.7}
\int_{\mfr{N}(\mbb{R}_+)}\e^{-\langle \nu, f\rangle}P_t(\sigma, d\nu)
 =
\exp\bigg\{-\<\sigma, u_tf\>-\int_0^t\psi(u_sf) d s\bigg\},
 \eeqlb
where $u_tf(x)$ is the unique solution to \eqref{2.3}. Such a process is characterized by the properties  (2.A), (2.B), (2.C) given in Section~2 and the following:
 \benumerate

\item[(3.A)] The immigrants come according to a Poisson random measure on $(0,\infty)\times \mfr{N}(\mbb{R}_+)^\circ$ with intensity $d sL(d\nu)$.

 \eenumerate

Let $M(dt,dy,dz,du,dv)$ be a Poisson random measure given as in Section~2 and let $N(dt,d\nu)$ be an $(\mcr{F}_t)$-Poisson random measure on $(0, \infty)\times \mfr{N}(\mbb{R}_+)^\circ$ with intensity $dtL(d\nu)$. We assume that the two random measures are independent with each other. Consider the following stochastic integral equation:
 \beqlb\label{2.11a}
Y_t(x) \ar=\ar Y_0\circ\theta_{-t}(x)+ \int_0^t\int_0^{\<Y_{s-}, \alpha\>}\int_0^1\int_{\mbb{N}}\int_0^{p({_\alpha \tilde{Y}}_{s-}^{-1}(y), z)} \Big[z\zeta_0 \cr
 \ar\ar\qqquad\qqquad\qquad
-\, \zeta_{{_\alpha \tilde{Y}}_{s-}^{-1}(y)}\Big]\circ\theta_{s-t}(x) M(ds,du,dy,dz,dv)\cr
 \ar\ar
+ \int_0^t\int_{\mfr{N}(\mbb{R}_+)^\circ} \nu\circ\theta_{s-t}(x)N(ds,d\nu),
 \eeqlb
where ${_\alpha \tilde{Y}}_{s-}^{-1}(y)= {_\alpha Y}_{s-}^{-1}( \<Y_{s-}, \alpha\>y).$ Here the first two terms on the right-hand side are as explained for \eqref{2.4added01}. The last term represents the immigration. A pathwise unique solution to \eqref{2.11a} is constructed as follows. Let $\sigma_0 = 0$. Given $\sigma_{k-1},$ we define
 \beqnn
\sigma_k = \sigma_{k-1} + \inf\{t>0: N((\sigma_{k-1},\sigma_{k-1}+t]\times N(\mbb{R}_+^\circ))> 0\},
 \eeqnn
and
 \beqnn
Y_t(x) = X_{k,t-\sigma_{k-1}}(x), \qquad \sigma_{k-1} \le t < \sigma_k,
 \eeqnn
where $\{X_{k,t}(x): t\ge 0\}$ is the pathwise unique solution to the following equation:
 \beqnn
X_t(x) \ar=\ar Y_{\sigma_{k-1}}\circ\theta_{-t}(x) + \int_0^t\int_0^{\<X_{s-},\alpha\>}\int_0^1\int_{\mbb{N}} \int_0^{p({_\alpha \tilde{X}}_{s-}^{-1}(y), z)} \Big[z \zeta_0\circ\theta_{s-t}(x) \cr
 \ar\ar\qqquad\qqquad\qquad
-\, \zeta_{{_\alpha \tilde{X}}_{s-}^{-1}(y)}\circ\theta_{s-t}(x)\Big] M(\sigma_{k-1}+ds,du,dy,dz,dv).
 \eeqnn
 Then we define
 \beqnn
Y_{\sigma_k}(x) = Y_{\sigma_k-}(x) + \int_{\{\sigma_k\}} \int_{\mfr{N}(\mbb{R}_+)^\circ} \nu(x) N(d s, d\nu).
 \eeqnn
It is easy to see that $\lim_{k\to \infty}\sigma_k = \infty$ and $\{Y_t: t \ge 0\}$ is the pathwise unique solution to \eqref{2.11a}.

The solution to \eqref{2.11a} determines a measure-valued strong Markov process $\{Y_t: t \ge 0\}$ with state space $\mfr{N}(\mbb{R}_+)$. We here omit the proofs of some of the following results since the arguments are similar to those for the corresponding results in Section~2.

\begin{prop}\label{bp3.1}
For any $f \in B(\mbb{R}_+)^+$ we have
 \beqlb\label{b3.18}
\langle Y_t, f\rangle \ar=\ar \langle Y_0, f\circ \theta_t\rangle + \int_0^t\int_0^{\<Y_{s-}, \alpha\>}\int_{[0, 1]}\int_{\mbb{N}}\int_0^{p({_\alpha \tilde{Y}}_{s-}^{-1}(y), z)}\Big[z f\circ\theta_{t - s}(0)
\cr
 \ar\ar\qqquad\qqquad\qquad
-\, f\circ\theta_{t - s}({_\alpha \tilde{Y}}_{s-}^{-1}(y))\Big] M(ds,du,dy,dz,dv) \cr
\ar\ar + \int_0^t\int_{\mfr{N}(\mbb{R}_+)^\circ} \<\nu, f\circ\theta_{t - s}\> N(d s, d\nu).
 \eeqlb
\end{prop}

\begin{prop}\label{bp3.1zz}
For any $t \ge 0$ and $f \in C^1(\mbb{R}_+)^+$ we have
 \beqlb\label{Yf}
\langle Y_t, f\rangle \ar=\ar \langle Y_0, f\rangle + \int_0^t \<Y_{s-}, f'\>d s + \int_0^t\int_0^{\<Y_{s-}, \alpha\>}\int_{[0, 1]}\int_{\mbb{N}}\int_0^{p({_\alpha \tilde{Y}}_{s-}^{-1}(y), z)}\Big[z f(0)
\cr
 \ar\ar\qqquad\qqquad\qquad
-\, f({_\alpha \tilde{Y}}_{s-}^{-1}(y))\Big] M(ds,du,dy,dz,dv) \cr
\ar\ar + \int_0^t\int_{\mfr{N}(\mbb{R}_+)^\circ} \<\nu, f\> N(d s, d\nu).
 \eeqlb
\end{prop}



\begin{prop}\label{bp3.3}
For any $f, G\in C^1(\mbb{R}_+)$ let $G_f(\mu)= G(\<\mu,f\>)$ and let
 \beqlb\label{b3.20}
LG_f(\mu)\ar=\ar \<\mu,f'\>G'(\<\mu,f\>) - \int_{\mbb{R}_+} \alpha(y) \sum_{z \in \mbb{N}}p(y,z)\Big[G(\<\mu, f\>) \ccr
 \ar\ar\qqquad\qqquad
-\, G\big(\<\mu, f\> + zf(0) - f(y)\big)\Big] \mu(d y) \cr
 \ar\ar\qquad
+ \int_{\mfr{N}(\mbb{R}_+)^\circ}\Big[G(\<\mu, f\> + \<\nu, f\>) - G(\<\mu, f\>)\Big]L(d\nu).
 \eeqlb
 Then we have
 \beqnn
 G(\<Y_t, f\>) = G(\<Y_0, f\>) + \int_0^t LG_f(Y_s)d s + \mbox{mart.}
 \eeqnn
 \end{prop}

\proof Let $\tilde{M}$ and $\tilde{N}$ denote the compensated measure of $M$ and $N$, respectively. By \eqref{Yf} and It\^{o}'s formula we have
{\small \beqnn
G(\<Y_t,f\>) \ar=\ar G(\<Y_0,f\>) + \int_0^t G'(\<Y_s, f\>)\<Y_s, f'\> d s \cr
\ar\ar + \int_0^t\int_0^{\<Y_{s-}, \alpha\>}\int_0^1\int_{\mbb{N}}\int_0^{p({_\alpha \tilde{Y}}_{s-}^{-1}(y), z)}\Big[G\big(\<Y_{s-},f\> + zf(0) - f({_\alpha \tilde{Y}}_{s-}^{-1}(y))\big)\cr
 \ar\ar\qqquad\qqquad\qquad
-\,G(\<Y_{s-},f\>) \Big]M(ds,du,dy,dz,dv)\cr
 \ar\ar+\int_0^t\int_{\mfr{N}(\mbb{R}_+)^\circ}\Big[G\big(\<Y_{s-},f\> + \<\nu, f\>\big) - G\big(\<Y_{s-},f\>\big)\Big]N(d s, d\nu)\cr
\ar=\ar G(\<Y_0,f\>) + \int_0^t G'(\<Y_s, f\>)\<Y_s, f'\> d s + N_t^G(f) \cr
 \ar\ar
+ \int_0^t \<Y_{s-},\alpha\>d s\int_{\mbb{R}_+} \alpha(y)\sum_{z\in \mbb{N}}p(y,z) \Big[G\big(\<Y_{s-},f\> + zf(0) - f(y)\big) \cr
 \cr
 \ar\ar\qqquad\qqquad\qquad
-\, G(\<Y_{s-},f\>) \Big]Y_{s-}(d y) \cr
\ar\ar + \int_0^t d t\int_{\mfr{N}(\mbb{R}_+)^\circ}[G\big(\<Y_{s-},f\> + \<\nu, f\>\big) - G(\<Y_{s-},f\>)]L(d\nu),
 \eeqnn}
where
\beqnn
N_t^G(f)\ar=\ar \int_0^t\int_0^{\<Y_{s-}, \alpha\>}\int_0^1\int_{\mbb{N}}\int_0^{p({_\alpha \tilde{Y}}_{s-}^{-1}(y), z)}\Big[G(\<Y_{s-},f\> + zf(0) - f({_\alpha \tilde{Y}}_{s-}^{-1}(y)))\cr
 \ar\ar\qqquad\qqquad\qquad
 -G(\<Y_{s-},f\>) \Big]\tilde{M}(ds,du,dy,dz,dv)\cr
 \ar\ar+\int_0^t\int_{\mfr{N}(\mbb{R}_+)^\circ}\Big[G\big(\<Y_{s-},f\> + \<\nu, f\>\big) - G\big(\<Y_{s-},f\>\big)\Big]\tilde{N}(d s, d\nu).
 \eeqnn
 By a first moment estimate one can check that $\{N_t^G(f): t\ge 0\}$ is a martingale. \qed

\begin{thm}\label{t2.2}
The measure-valued process $\{Y_t: t \ge 0\}$ defined by \eqref{2.11a} is an $(\alpha, g, \psi)$- age-structured branching process with immigration.
\end{thm}

\proof
Let $G \in C^1(\mbb{R}_+)$ and let $t\to f_t$ be a mapping from $[0, T]$ to $C^1(\mbb{R}_+)^+$ such that $t\to f_t$ is continuously differentiable and $t\to f_t'$ is continuous by the supremum norm. Using Proposition~\ref{bp3.3} one can show as in the proof of Theorem~\ref{t2.2a} that
 \beqnn
G(\<Y_t,f_t\>) \ar=\ar G(\<Y_0,f_0\>) + \int_0^t \Big\{ G'(\<Y_s,f_s\>)\<Y_s,f_s'\> + \int_{\mbb{R}_+}\alpha(y)Y_s(d y) \cr
\ar\ar \sum_{z \in \mbb{N}}p(y,z)\Big[G(\<Y_s,f_s\> + zf_s(0) - f_s(y)) - G(\<Y_s,f_s\>)\Big]\Big\}d s \cr
\ar\ar + \int_0^t d s \int_{\mfr{N}(\mbb{R}_+)^\circ}[G(\<Y_{s},f_s\> + \<\nu,f_s\>) - G(\<Y_{s-},f_s\>)]L(d\nu) \cr
\ar\ar + \int_0^t G'(\<Y_s,f_s\>)\<Y_s,f_s'\> d s + N(t),
 \eeqnn
where $\{N(t)\}$ is a martingale. For any $f\in C^1(\mbb{R}_+)^+$ we may apply the above relation to $G(z) = e^{-z}$ and
 $$
f_t = u_{T -t}f + \int_0^{T-t}d s \int_{\mfr{N}(\mbb{R}_+)^\circ}[1 - e^{-\<\nu, u_sf\>}] L(d\nu).
 $$
In this case, one can see that
 \beqnn
t\mapsto \exp\Big\{-\<Y_t, u_{T-t}f\> - \int_0^{T-t}\psi(u_sf)d s\Big\}
 \eeqnn
is a martingale. Then $\{Y_t: t\ge 0\}$ has transition semigroup $(Q_t)_{t\ge 0}$ defined by \eqref{1.7}. That ends the proof. \qed

A necessary and sufficient condition for the ergodicity of the $(\alpha, g, \psi)$-age-structured branching process with immigration is given by the following theorem:

\begin{thm}\label{t2.1}
Suppose that $c_0= \sup_{y\ge 0}\alpha(y)[g'(y,1-)-1]< 0$. Then $P_t(\sigma, \cdot)$ converges to a probability measure $\eta$ on $\mfr{N}(\mbb{R}_+)$ as $t\to \infty$ for every $\sigma \in \mfr{N}(\mbb{R}_+)$ if and only if
 \beqlb\label{t2.1ss}
\int_{\mfr{N}(\mbb{R}_+)} 1_{\{\<\nu,1\>\ge 1\}}\log \<\nu,1\> L(d\nu) < \infty.
 \eeqlb
In this case, the Laplace transform of $\eta$ is given by
 \beqlb\label{t2.1tt}
\int_{\mfr{N}(\mbb{R}_+)} e^{-\<\nu, f\>} \eta(d\nu) = \exp\bigg\{- \int_0^\infty \psi(u_sf)d s\bigg\}.
 \eeqlb
\end{thm}

\proof Let $f_*= \inf_{x\ge 0} f(x)$ for $f\in B(\mbb{R}_+)^+$. By Jensen's inequality one can check that
\beqnn
e^{-u_tf(x)} = \mbf{E}e^{-\<X_t^{\delta_x}, f\>} \ge e^{-\mbf{E} \<X_t^{\delta_x}, f\>} = e^{-\pi_tf(x)},
\eeqnn
which implies $u_tf(x) \le \pi_tf(x)$ for any $x \ge 0$ and $f \in B(\mbb{R}_+)^+.$ By Proposition~\ref{t2.1b} and Corollary~\ref{t2.1c}, we have
 \beqlb\label{2.5d}
(1-\e^{-f_*})e^{-c_1t}\le \|\pi_tf\|\le   e^{c_0t}\|f\|.
 \eeqlb
For any $a,c>0$ it is elementary to see that
 \beqlb\label{2.5f}
\ar\ar\int_0^\infty d s\int_{\mfr{N}(\mbb{R}_+)^\circ} (1-\e^{-a\<\nu,1\>\e^{-cs}})L(d\nu) \cr
 \ar\ar\qqquad
= \int_{\mfr{N}(\mbb{R}_+)^\circ}L(d\nu) \int_0^\infty (1-\e^{-a\<\nu,1\>\e^{-cs}})d s \cr
 \ar\ar\qqquad
= c^{-1}\int_{\mfr{N}(\mbb{R}_+)^\circ}L(d\nu) \int_0^{c\<\nu,1\>} (1-\e^{-z})z^{-1}d z.
 \eeqlb
Clearly, the above quantity is finite if and only if \eqref{t2.1ss} holds. Suppose that \eqref{t2.1ss} holds. Using \eqref{2.5d} and dominated convergence we see the functional
 \beqlb\label{2.5e}
\int_0^\infty \psi(u_sf)d s
 =
\int_0^\infty d s\int_{\mfr{N}(\mbb{R}_+)^\circ} (1-e^{-\<\nu,u_sf\>})L(d\nu)
 \eeqlb
is continuous on $B(\mbb{R}_+)^+$ with respect to the pointwise convergence. By \eqref{1.7} we infer that $P_t(\sigma,\cdot)$ converges to a probability measure $\eta$ defined by \eqref{t2.1tt}; see, e.g., Theorem~1.20 in Li (2011). Conversely, suppose that $P_t(\sigma, \cdot)$ converges to a probability measure $\eta$ on $\mfr{N}(\mbb{R}_+)$ as $t\to \infty$ for every $\sigma \in \mfr{N}(\mbb{R}_+)$. By \eqref{1.7} we see that $\eta$ has Laplace functional given by \eqref{t2.1tt}. Then \eqref{2.5e} is finite for every $f\in B(\mbb{R}_+)^+$. Taking $f\equiv 1$ and using \eqref{2.5d} we see the quantity in \eqref{2.5f} is finite. Then we have \eqref{t2.1ss}. \qed

\bigskip

{\bf Acknowledgement} The research of Li was supported by NSFC (No.\,11531001 and No.\,11626245). Ji's research was supported by China Postdoctoral Science Foundation (No.\,2020M681994). We are very thankful to Professor Jie Xiong for helpful discussions.

\end{document}